\documentclass{amsart}

\usepackage{amsmath}
\usepackage{amscd}
\usepackage{amssymb}
 \usepackage{color}

\input xy
\xyoption{all}

\newcommand{\pd}{\partial}
\newcommand{\bC}{{\mathbb C}}

\newcommand{\bZ}{{\mathbb Z}}

 \newcommand{\cD}{{\mathcal D}}

\newcommand{\cL}{{\mathcal L}}

\newtheorem{theorem/definition}{Theorem/Definition}[section]

\newtheorem{Conjecture}{Conjecture}

\theoremstyle{remark}

\theoremstyle{definition}

\newcommand{\bml}{\begin{multline}}
\newcommand{\eml}{\end{multline}}

\newcommand{\bag}{\begin{align}}
\newcommand{\egn}{ \end{align}}

\newcommand{\be}{\begin{equation}}
\newcommand{\ee}{\end{equation}}
\newcommand{\bea}{\begin{eqnarray}}
\newcommand{\eea}{\end{eqnarray}}
\newcommand{\ben}{\begin{eqnarray*}}
\newcommand{\een}{\end{eqnarray*}}
\newcommand{\bet}{\begin{equation}
\begin{split}}
\newcommand{\eet}{\end{split}
\end{equation}}

\definecolor{yellow}{rgb}{1,1,0}
\definecolor{orange}{rgb}{1,.7,0}
\definecolor{red}{rgb}{1,0,0} \definecolor{blue}{rgb}{0,0,1}
\definecolor{white}{rgb}{1,1,1}

\definecolor{A}{rgb}{.75,1,.75}

\begin{document}

\title
{Integrality Properties of Open-Closed Mirror Maps}
\author{Jian Zhou}
\address{Department of Mathematical Sciences\\Tsinghua University\\Beijing, 100084, China}
\email{jzhou@math.tsinghua.edu.cn}

\begin{abstract}
We propose a conjecture on integrality property of the open-closed mirror maps
of compact Calabi-Yau manifolds.
Some examples are presented.
\end{abstract}


\maketitle

In this paper we will propose a conjecture on integrality properties
of open-closed mirror maps for compact Calabi-Yau manifolds,
inspired by the Lian-Yau integrality for closed mirror map for compact Calabi-Yau manifold
\cite{Lia-Yau, Lia-Yau2, Zud, Kra-Riv, Kra-Riv2, Kra-Riv3, Del},
and its recent extension to the local Calabi-Yau case \cite{Zho1}.
The integrality of open-closed mirror map in the noncompact Calabi-Yau $3$-fold case
was observed in \cite{Zho1},
in this work we propose to extend it to the compact Calabi-Yau case and more noncompact examples, not necessarily three dimensional,
based on a related work \cite{Zho2}.

\section{Integrality Conjecture of Open-Closed Mirror Maps}

\subsection{The conjecture}

Suppose that we have a system of Picard-Fuchs equations in variables $x_0, x_1, \dots, x_N$
that comes from some charge vectors
describing some brane geometry in some Calabi-Yau geometry,
in the large volume phase.
Then it is expected that the system has a holomorphic solution
$g_0(x_0, x_1, \dots, x_N)$ such that $g_0(0, \dots, 0) = 1$
and $N$ logarithmic solutions of the form
$$g_1^{(i)}(x_0, \dots, x_N) = \log x_i \cdot g_0(x_0, \dots, x_N)
+ h_1^{(i)}(x_0, \dots, x_N), \qquad i =0, 1, \dots, N,$$
where $h_1^{(i)}$ are holomorphic functions such that
$h^{(i)}_1(0, \dots, 0) = 0$.
The open-closed mirror map is defined by:
\be
q_i = \exp(g_1^{(i)}/g_0) = x_i \exp (h_1^{(i)}/g_0), \qquad i=0, 1, \dots, N.
\ee
(For some examples,
see e.g. \cite{Ler-May}.)
On can find the inverse map by Lagrange-Good inversion formula \cite{Good} as in \cite{Zho1}.

Our conjecture is

\begin{Conjecture}
The Taylor series of $q_0, q_1, \dots, q_N$ in $x_0, x_1, \dots, x_N$ are integers,
and vice versa.
\end{Conjecture}

We will present some examples for which we verify this conjecture.
We will focus on some cases with one closed string moduli parameter.
Extensions to the case of multiple closed string moduli parameters
will be left to subsequent work.

\subsection{Some Picard-Fuchs systems related to open-closed mirror map}

Let $k_1,  \dots, k_n$ be positive integers such that
\be \label{eqn:Sum=1}
\frac{1}{k_1} + \cdots + \frac{1}{k_n} = 1.
\ee
See \cite{Zho2} for a list of solutions for $n=2,3,4,5$
under the assumption that $k_1 \leq  \dots \leq k_n$.
Here we do not require this assumption.
Let $k=l.c.m.(k_1, \dots, k_n)$ and $w_i = \frac{k}{k_i}$.
Then we have
\be
k_1 w_1 = \cdots = k_n w_n = w_1 + \cdots + w_n = k.
\ee
Consider charge vectors formed from such a solution:
\be \label{eqn:Charge0}
\begin{pmatrix}
\tilde{l}^{(1)} \\ \tilde{l}^{(0)}
\end{pmatrix}
= \biggl( \begin{array}{c|cccc|cc}
-k & w_1 & w_2 & \cdots & w_n & 0 & 0 \\
-1 &  1 & 0 & \cdots & 0 & -1 & 1
\end{array} \biggr).
\ee
We will transform them to
\be
\begin{pmatrix}
l^{(1)} \\ l^{(0)}
\end{pmatrix}
= \biggl( \begin{array}{c|cccc|cc}
-k+w_1 & 0 & w_2 & \cdots & w_n & w_1 & -w_1 \\
-1 &  1 & 0 & \cdots & 0 & -1 & 1
\end{array} \biggr),
\ee
where $l^{(1)} = \tilde{l}^{(1)} - w_1 \tilde{l}^{(0)}$ and $l^{(0)} = \tilde{l}^{(0)}$.
We will consider the extended Picard-Fuchs system \cite{Ler-May} associated with three charge vectors
$l^{(1)}$, $l^{(0)}$ and $\tilde{l}_1 = l^{(1)} + w_1 l^{(0)}$:
\be \label{eqn:PF}
\cL_1 S(x_0,x_1) =0, \qquad \cL_0 S(x_0,x_1) = 0, \qquad \cL_1' S(x_0,x_1) = 0,
\ee
where
\bea
&& \cL_0 =  \big[ \theta_0 - x_0 (1 + \theta_0 + (k-w_1) \theta_1) \big] (\theta_0 - w_1 \theta_1), \\
&& \cL_1 = \prod_{i=2}^n \prod_{j=0}^{w_i-1} (w_i \theta_1-j)
\prod_{j=0}^{w_1-1} (-\theta_0+w_1\theta_1-j) \\
&& \qquad - x_1 \prod_{j=1}^{k-w_1} (\theta_0+(k-w_1)\theta_1+j)
\prod_{j=0}^{w_1-1} (-\theta_0+w_1\theta_1+j), \nonumber \\
&& \cL_1'
= \prod_{j=0}^{w_1-1} (\theta_0-j) \cdot \prod_{i=2}^n \prod_{j=0}^{w_i-1} (w_i \theta_1 -j)
- x_0^{w_1}x_1 \prod_{j=1}^{k} (\theta_0+(k-w_1)\theta_1+j),
\eea
where $\theta_i = x_i \frac{\pd}{\pd x_i}$, $i=0,1$.
When $n=4$,
some examples of this system correspond to the extended Picard-Fuchs system associated with
B-branes in compact Calabi-Yau $3$-folds in the large volume phase (see e.g. \cite{AHMM}).
See  \cite{Li-Lia-Yau} and its references for some recent work on extension of mirror symmetry to
the open string case.

\subsection{Derivation of the Picard-Fuchs operators}
Let us recall the well-known procedure of derivation of Picard-Fuchs operators from the charge vectors
for the convenience of the reader.
Given a vector $l = (l_0, l_1, \dots, l_{n+2}) \in \bZ^{n+3}$,
consider the affine complex space $\bC^{p+1}$ with linear coordinates $(a_0, \dots, a_{n+2})$,
define an operator
\be
\cD_l = \prod_{l_j > 0} \big( \frac{\pd}{\pd a_i} \big)^{l_j}
- \prod_{l_j < 0} \big( \frac{\pd}{\pd a_i} \big)^{-l_j}.
\ee
Introduce variables
\be
x_0 = \frac{a_1a_{n+2}}{a_0a_{n+1}}, \qquad
x_1 = (-1)^k \frac{a_2^{w_2} \cdots a_n^{w_n}a_{n+1}^{w_1}}{a_0^{k-w_1} a_{n+2}^{w_1}}.
\ee
Suppose that we have a function $\tilde{\Pi}(x_0,x_1)$ such that
\be
\Pi(a_0, \dots, a_{n+2}) = \frac{1}{a_0} \tilde{\Pi}(x_0,x_1)
\ee
satisfies the equations
\be
\cD_{l^{(i)}} \Pi(a_0, \dots, a_{n+2}) =0, \quad i =0, 1, \qquad
\cD_{\tilde{l}^{(1)}} \Pi(a_0, \dots, a_{n+2}) = 0.
\ee
For $i=0$,
\ben
&& \prod_{l_j^{(0)}>0} \big( \frac{\pd}{\pd a_i} \big)^{l^{(0)}_j} \Pi(a_0, \dots, a_{n+2})
= \frac{\pd}{\pd a_1} \frac{\pd}{\pd a_{n+2}} (\frac{1}{a_0} \tilde{\Pi}(x_0, x_1)) \\
& = & \frac{1}{a_0} \frac{\pd}{\pd a_1} (\frac{1}{x_0} \frac{\pd x_0}{\pd a_{n+2}} \theta_0 \tilde{\Pi}(x_0,x_1)
+ \frac{1}{x_1} \frac{\pd x_1}{\pd a_{n+2}} \theta_1 \tilde{\Pi}(x_0,x_1)) \\
& = & \frac{1}{a_0a_{n+2}} \frac{\pd}{\pd a_1} ( \theta_0 \tilde{\Pi}(x_0,x_1)
- w_1 \theta_1 \tilde{\Pi}(x_0,x_1)) \\
& = &  \frac{1}{a_0a_{n+2}}
(\frac{1}{x_0} \frac{\pd x_0}{\pd a_1} \theta_0
+ \frac{1}{x_1} \frac{\pd x_1}{\pd a_1} \theta_1) ( \theta_0 - w_1 \theta_1) \tilde{\Pi}(x_0,x_1)) \\
& = & \frac{1}{a_0a_1a_{n+2}} \theta_0 (\theta_0 - w_1 \theta_1) \tilde{\Pi}(x_0, x_1),
\een
and
\ben
&& \prod_{l_j < 0} \big( \frac{\pd}{\pd a_i} \big)^{-l^{(0)}_j}
\Pi(a_0, \dots, a_{n+2})
= \frac{\pd}{\pd a_0} \frac{\pd}{\pd a_{n+1}} (\frac{1}{a_0} \tilde{\Pi}(x_0,x_1)) \\
& = & \frac{\pd}{\pd a_0} (\frac{1}{a_0} \frac{1}{x_0}\frac{\pd x_0}{\pd a_{n+1}} \theta_0 \tilde{\Pi}(x_0,x_1)
+ \frac{1}{a_0} \frac{1}{x_1}\frac{\pd x_1}{\pd a_{n+1}} \theta_1 \tilde{\Pi}(x_0,x_1)) \\
& = & \frac{1}{a_{n+1}}\frac{\pd}{\pd a_0} (\frac{1}{a_0} \cdot
(-\theta_0 + w_1 \theta_1) \tilde{\Pi}(x_0,x_1)) \\
& = & \frac{1}{a_{n+1}} (-\frac{1}{a_0^2} + \frac{1}{a_0}  \frac{1}{x_0}\frac{\pd x_0}{\pd a_0} \theta_0
+ \frac{1}{a_0}  \frac{1}{x_1}\frac{\pd x_1}{\pd a_0} \theta_1)
(-\theta_0 + w_1 \theta_1) \tilde{\Pi}(x_0,x_1) \\
& = &  \frac{1}{a_0^2a_{n+1}} (-1 - \theta_0 - (k-w_1) \theta_1)
(-\theta_0 + w_1 \theta_1) \tilde{\Pi}(x_0,x_1).
\een
Therefore,
\ben
\theta_0 (\theta_0 - w_1 \theta_1) \tilde{\Pi}(x_0, x_1)
+ \frac{a_1a_{n+2}}{a_0a_{n+1}} (1 + \theta_0 + (k-w_1) \theta_1)
(-\theta_0 + w_1 \theta_1) \tilde{\Pi}(x_0,x_1) = 0.
\een
This is the equation $\cL_0 \tilde{\Pi}(x_0,x_1) = 0$.
The other two equations are derived in the same fashion.

\subsection{Holomorphic solutions}

In this subsection we will show that the extended Picard-Fuch system \eqref{eqn:PF}
has a holomorphic solution of the form
$$g_0 = 1 + \sum_{m_0+ m_1 > 0} a_{m_0,m_1} x_0^{m_0}x_1^{m_1}.$$
We will derive its explicit expression as follows.

The Picard-Fuchs equations for $g_0$ can be rewritten as follows:
\ben
&& \sum_{m_0+ m_1 > 0} m_0(m_0-w_1m_1) a_{m_0,m_1} x_0^{m_0}x_1^{m_1} \\
& = & x_0 \sum_{m_0+ m_1 > 0} (m_0-w_1m_1) (1+m_0+(k-w_1)m_1) a_{m_0,m_1} x_0^{m_0}x_1^{m_1}, \\
&& \sum_{m_0+ m_1 > 0} \prod_{j=0}^{w_1-1} (-m_0+w_1m_1+j) \prod_{i=2}^n \prod_{j=0}^{w_i-1}
(w_i m_1-j) a_{m_0,m_1} x_0^{m_0}x_1^{m_1} \\
& = &  x_1 \prod_{j=1}^{k-w_1} (m_0+(k-w_1)m_1+j)
\prod_{j=0}^{w_1-1} (-m_0+w_1m_1+j) a_{m_0,m_1} x_0^{m_0}x_1^{m_1}, \\
&& \sum_{m_0+ m_1 > 0} \prod_{j=0}^{w_1-1} (m_0-j) \cdot \prod_{i=2}^n \prod_{j=0}^{w_i-1} (w_i m_1 -j) a_{m_0,m_1} x_0^{m_0}x_1^{m_1} \\
& = & x_0^{w_1}x_1(k!+ \sum_{m_0+ m_1 > 0} \prod_{j=1}^{k} (m_0+(k-w_1)m_1+j) a_{m_0,m_1} x_0^{m_0}x_1^{m_1}).
\een
These equations give rise to some recursions which can be solved directly to give us
the following solution:
\be
g_0(x_0,x_1) = \sum_{m \geq 0} \frac{(km)!}{\prod_{i=1}^n (w_im)!} (x_0^{w_1}x_1)^m.
\ee
Clearly $g_0(x_0,x_1) \in \bZ[[x_0^{w_1}x_1]]$.

\subsection{Logarithmic solutions}

In this subsection we will show that the extended Picard-Fuch system \eqref{eqn:PF}
has two logarithmic solutions of the form
$$g_1^{(i)} = g_0 \log x_i + \sum_{m_0+ m_1 > 0} a^{(i)}_{m_0,m_1} x_0^{m_0}x_1^{m_1}$$ ($i=0,1$)
to the Picard-Fuchs system \eqref{eqn:PF}.
We will derive their explicit expressions as follows.
From $\cL_1'g_1^{(1)} =0$ we get
\ben
&& \sum_{m_0+ m_1 > 0}  \prod_{j=0}^{w_1-1} (m_0-j) \cdot \prod_{i=2}^n \prod_{j=0}^{w_i-1} (w_i m_1 -j) a^{(1)}_{m_0,m_1} x_0^{m_0}x_1^{m_1} \\
& + & \sum_{m=1}^\infty \frac{(km)!}{\prod_{i=1}^n (w_i(m-1))!}
\sum_{i=2}^n \sum_{j=0}^{w_i-1} \frac{1}{w_im-j} (x_0^{w_1} x_1)^m\\
& = & x_0^{w_1}x_1 \sum_{m_0+ m_1 > 0} \prod_{j=1}^{k} (m_0+(k-w_1)m_1+j) a^{(1)}_{m_0,m_1} x_0^{m_0}x_1^{m_1} \\
& + & x_0^{w_1}x_1\sum_{m=0}^\infty \frac{(k(m+1))!}{\prod_{i=1}^n (w_im)!}
\sum_{j=1}^k \frac{k-w_1}{km+j} (x_0^{w_1}x_1)^m.
\een
Compare the coefficients of $x_0^{m_0}x_1^{m_1}$ on both sides.
When $w_0 = w_1 m$ and $w_2 = m$,
\ben
&& a^{(1)}_{w_1m, m}
+ \frac{(km)!}{\prod_{i=1}^n (w_im)!} \sum_{i=2}^n \sum_{j=0}^{w_i-1} \frac{1}{w_im-j}  \\
& = &  \frac{\prod_{j=0}^{k-1} (km-j)}{\prod_{i=1}^n \prod_{j=0}^{w_i-1} (w_i m -j)} a^{(1)}_{w_1(m-1),m-1}
+  \frac{(km)!}{\prod_{i=1}^n (w_im))!}\sum_{j=0}^{k-1} \frac{k-w_1}{km-j}.
\een
Hence we have
\be
a^{(1)}_{w_1m,m} =  \frac{(km)!}{\prod_{i=1}^n (w_im)!}
(\sum_{j=1}^{km} \frac{k-w_1}{j} - \sum_{i=2}^n \sum_{j=1}^{w_im} \frac{1}{j}).
\ee
For $m_1 \neq w_1m_0$,
we have
\ben
&& \prod_{j=0}^{w_1-1} (m_0-j) \cdot \prod_{i=2}^n \prod_{j=0}^{w_i-1} (w_i m_1 -j) a^{(1)}_{m_0,m_1}  \\
& = & \prod_{j=0}^{k-1} (m_0+(k-w_1)m_1-j) a^{(1)}_{m_0-w_1,m_1-1},
\een
hence when $m_0 > w_1m_1$ we have
\be \label{eqn:>}
a^{(1)}_{m_0,m_1} = \frac{(m_0+(k-w_1)m_1)!}{m_0! \prod_{i=2}^n (w_im_1)!}
a^{(1)}_{m_0-w_1m_1, 0}.
\ee
We will determine $a^{(1)}_{m,0}$ from another equation.

From $\cL_0 g_1^{(1)} = 0$ we get
\ben
&& \sum_{m_0+ m_1 > 0} m_0(m_0-w_1m_1) a^{(1)}_{m_0,m_1} x_0^{m_0}x_1^{m_1}
- w_1  \sum_{m \geq 1} \frac{(km)!w_1m }{\prod_{i=1}^n (w_im)!} (x_0^{w_1}x_1)^m \\
& = & x_0 \sum_{m_0+ m_1 > 0} (m_0-w_1m_1) (1+m_0+(k-w_1)m_1) a^{(1)}_{m_0,m_1} x_0^{m_0}x_1^{m_1} \\
& - & w_1 x_0 \sum_{m \geq 0} \frac{(km+1)!}{\prod_{i=1}^n (w_im)!} (x_0^{w_1}x_1)^m.
\een
Compare the coefficients of $x_0^{m_0}x_1^{m_1}$ on both sides.
For $m_0=m>0$ and $m_1 = 0$,
\ben
&& a_{1,0}^{(1)} = - w_1, \qquad
 m^2 a^{(1)}_{m,0} = m (m-1) a^{(1)}_{m-1,0},
\een
and so
\ben
a^{(1)}_{m,0} = - \frac{1}{m}w_1.
\een
Hence by combining with \eqref{eqn:>} for $m_0 > w_1m_1$,
\be
a^{(1)}_{m_0,m_1} = - w_1 \frac{(m_0+(k-w_1)m_1)!}{m_0! \prod_{i=2}^n (w_im_1)!} \frac{1}{m_0 - w_1m_1}.
\ee
When $m_0 \neq 0, w_1m_1$ or $w_1m_1+1$,
\be
a^{(1)}_{m_0,m_1}
=  \frac{(m_0-1-w_1m_1) (m_0+(k-w_1)m_1)}{m_0 (m_0-w_1m_1)} a^{(1)}_{m_0-1,m_1},
\ee
therefore,
for $m_0 < w_1m_1$,
we have
\be \label{eqn:<}
a^{(1)}_{m_0,m_1}
= \frac{(m_0+(k-w_1)m_1)!}{m_0!((k-w_1)m_1)!} \frac{w_1m_1}{w_1m_1-m_0} a^{(1)}_{0,m_1}.
\ee
We will determine $a^{(1)}_{0,m_1}$ from another equation.

From $\cL_1 g_1^{(1)} = 0$ we get:
\ben
&& \sum_{m_0+ m_1 > 0} \prod_{j=0}^{w_1-1} (-m_0+w_1m_1-j) \prod_{i=2}^n \prod_{j=0}^{w_i-1}
(w_i m_1-j) a^{(1)}_{m_0,m_1} x_0^{m_0}x_1^{m_1} \\
& + &  \sum_{m=1}^\infty
\frac{(-1)^{w_1-1} w_1! (km)!}{(w_1m)! \prod_{i=2}^n (w_i(m-1))!} (x_0^{w_1}x_1)^m \\
& = &   x_1 \sum_{m_0+ m_1 > 0} \prod_{j=1}^{k-w_1} (m_0+(k-w_1)m_1+j)
\prod_{j=0}^{w_1-1} (-m_0+w_1m_1+j) a^{(1)}_{m_0,m_1} x_0^{m_0}x_1^{m_1}\\
& + & x_1 \sum_{m=1}^\infty
\frac{w_1! (km+k-w_1)!}{\prod_{i=1}^n (w_im)!} (x_0^{w_1}x_1)^m.
\een
By comparing the coefficients of $x_1^{m_1}$ we get:
\ben
&& \prod_{i=1}^n w_i! \cdot a^{(1)}_{0,1} =  w_1! (k-w_1)!, \\
&& \prod_{i=1}^n \prod_{j=0}^{w_i-1}
(w_i m_1-j) a^{(1)}_{0,m_1}  =  \prod_{j=0}^{k-w_1-1} ((k-w_1)m_1-j)
\prod_{j=1}^{w_1} (w_1m_1-j) a^{(1)}_{0,m_1-1},
\een
hence
\be
a^{(1)}_{0,m_1} =  \frac{((k-w_1)m_1)!}{\prod_{i=2}^n (w_im_1)!} \frac{1}{m_1}.
\ee
Combining with \eqref{eqn:<} we get:
\be
a^{(1)}_{m_0,m_1} =- w_1 \frac{(m_0+(k-w_1)m_1)!}{m_0!\prod_{i=2}^n (w_im)!} \frac{1}{m_0-w_1m_1}.
\ee

Therefore,
we get the following solution:
\bea
g_1^{(1)}(x_0,x_1) & = & \log x_1 \cdot \sum_{m \geq 0} \frac{(km)!}{\prod_{i=1}^n (w_im)!} (x_0^{w_1}x_1)^m \\
& + & \sum_{m=1}^\infty \frac{(km)!}{\prod_{i=1}^n (w_im)!}
(\sum_{j=1}^{km} \frac{k-w_1}{j} - \sum_{i=2}^n \sum_{j=1}^{w_im} \frac{1}{j})(x_0^{w_1}x_1)^m  \nonumber \\
& - & \sum_{m_0+m_1>0, m_0 \neq w_1m_1}
w_1 \frac{(m_0+(k-w_1)m_1)!}{m_0!\prod_{i=2}^n (w_im)!} \frac{x_0^{m_0} x_1^{m_1}}{m_0-w_1m_1}. \nonumber
\eea
In the same fashion one can get the following solution:
\bea
g_1^{(0)}(x_0,x_1) & = & \log x_0 \cdot \sum_{m \geq 0} \frac{(km)!}{\prod_{i=1}^n (w_im)!} (x_0^{w_1}x_1)^m \\
& + & \sum_{m=1}^\infty \frac{(km)!}{\prod_{i=1}^n (w_im)!}
(\sum_{j=1}^{km} \frac{1}{j} - \sum_{j=1}^{m} \frac{1}{w_1j}) (x_0^{w_1}x_1)^m \nonumber \\
& + & \sum_{m_0+m_1>0, m_0 \neq w_1m_1}
\frac{(m_0+(k-w_1)m_1)!}{m_0!\prod_{i=2}^n (w_im)!} \frac{x_0^{m_0} x_1^{m_1}}{m_0-w_1m_1}. \nonumber
\eea

\subsection{Open-closed mirror map and its integrality}

The open-closed mirror map is given by
\be
q_0 = \exp (g_1^{(0)}/g_0), \qquad
q_1 = \exp (g_1^{(1)} /g_0).
\ee
As a special case of Conjecture 1,
one should have
\be
q_0, q_1 \in \bZ[[x_0,x_1]].
\ee
Notice that in both cases we have some extra terms,
so the integrality property of open-closed mirror map is different from
the Lian-Yau integrality of closed mirror map,
because it involves the integrality of the following series:
\be
\exp (\sum_{m_0 \neq w_1m_1}
\frac{(m_0+(k-w_1)m_1)!}{m_0!\prod_{i=2}^n (w_im)!} \frac{x_0^{m_0} x_1^{m_1}}{m_0-w_1m_1}/
\sum_{m \geq 0} \frac{(km)!}{\prod_{i=1}^n (w_im)!} (x_0^{w_1}x_1)^m).
\ee

Applying the Lagrange-Good inversion formula as in \cite{Zho1},
one can get the explicit expressions for the inverse open-closed mirror map.
Write
\be
q_i = \sum_{m_0+m_1>0} A_{m_0,m_1}^{(i)} x_0^{m_0}x_1^{m_1},
\ee
then $A_{m_0,m_1}^{(0)}$ and $A_{m_0,m_1}^{(1)}$ are the coefficients of $x_0^{m_0-1}x_1^{m_1}$
and $x_0^{m_0}x_1^{m_1-1}$ of
\be
(1+\theta_0 h_0 + \theta_1 h_1) \cdot \exp (-m_0h_0-m_1h_1),
\ee
where $h_i = g_1^{(i)}/g_0 - \log x_i$.

\subsection{The Picard-Fuchs system in a different phase}
One can also consider the Picard-Fuchs system associated with the charge vectors \eqref{eqn:Charge0}.
The associated variables are:
\be
\tilde{x}_0 = \frac{a_1a_{n+2}}{a_0a_{n+1}}, \qquad \tilde{x}_1 = (-1)^k \frac{a_1^{w_1} \cdots a_n^{w_n}}{a_0^k}.
\ee
I.e.,
\be \label{eqn:Change}
\tilde{x}_0 = x_0, \qquad \tilde{x}_1 = x_0^{w_1} x_1.
\ee
Let $\tilde{\theta}_i = \tilde{x}_i \frac{\pd}{\pd \tilde{x}_i}$.
Then we have:
\begin{align}
\theta_0 & = \tilde{\theta}_0 + w_1 \tilde{\theta}_1, &
 \theta_1 & = \tilde{\theta}_1, \\
\tilde{\theta}_0 & = \theta_0 - w_1 \theta_1, & \tilde{\theta}_1 & = \theta_1.
\end{align}
Under such change of variables,
the Picard-Fuchs operators become
\bea
&& \cL_0 =  \big[ \tilde{\theta}_0  + w_1 \tilde{\theta}_1- x_0 (1 + \tilde{\theta}_0 + k \tilde{\theta}_1) \big] \tilde{\theta}_0, \\
&& \cL_1 = \prod_{i=2}^n \prod_{j=0}^{w_i-1} (w_i \tilde{\theta}_1-j)
\prod_{j=0}^{w_1-1} (-\tilde{\theta}_0-j) \\
&& \qquad - \tilde{x}_0^{-w_1} \tilde{x}_1 \prod_{j=1}^{k-w_1} (\tilde{\theta}_0+k\tilde{\theta}_1+j)
\prod_{j=0}^{w_1-1} (-\tilde{\theta}_0+j), \nonumber \\
&& \cL_1'
= \prod_{j=0}^{w_1-1} (\tilde{\theta}_0+w_1\tilde{\theta}_1-j) \cdot
\prod_{i=2}^n \prod_{j=0}^{w_i-1} (w_i \tilde{\theta}_1 -j)
- \tilde{x}_1 \prod_{j=1}^{k} (\tilde{\theta}_0+k\tilde{\theta}_1+j).
\eea
Consider the system given by
\be
\cL_0 S(\tilde{x}_0, \tilde{x}_1) = 0, \qquad \cL_1' S(\tilde{x}_0, \tilde{x}_1) = 0.
\ee
It has a solution holomorphic at $\tilde{x}_0 = \tilde{x}_1 = 0$:
\be
\tilde{g}_0(\tilde{x}_0, \tilde{x}_1)
= 1 + \sum_{m=1}^\infty \frac{(km)!}{\prod_{i=1}^n (w_im)!} \tilde{x}_1^m,
\ee
and a logarithmic solution:
\be
\tilde{g}_1^{(1)}(\tilde{x}_0, \tilde{x}_1)
= \log \tilde{x}_1 \cdot \tilde{g}_0(\tilde{x}_0, \tilde{x}_1)
+ \sum_{m=1}^\infty \frac{(km)!}{\prod_{i=1}^n (w_im)!}
(\sum_{j=1}^{km} \frac{k}{j} - \sum_{i=1}^n \sum_{j=1}^{w_im} \frac{1}{j}) \tilde{x}_1^m.
\ee
They are not obtained from $g_0(x_0,x_1)$ and $g_1^{(1)}(x_0,x_1)$ by change of variables \eqref{eqn:Change},
and they do not satisfy the equation
$\cL_1 S(\tilde{x}_0, \tilde{x}_1) = 0$.
Furthermore,
one can easily see that there is no solution of the form
$\tilde{g}_1^{(0)}
= \log \tilde{x}_0 \cdot \tilde{g}_0(\tilde{x}_0, \tilde{x}_1)
+ \tilde{h}_1^{(1)}(\tilde{x}_0, \tilde{x}_1)$,
with $\tilde{h}_1^{(1)}$ holomorphic in $\tilde{x}_0$ and $\tilde{x}_1$.
Indeed,
if we write $\tilde{h}_1^{(0)}(\tilde{x}_0, \tilde{x}_1)
= \sum_{m_0+m_1> 0} \tilde{a}_{m_0,m_1}^{(0)} \tilde{x}_0^{m_0} \tilde{x}_1^{m_1}$,
then from $\cL_1' \tilde{h}_1^{(0)} = 0$ we get:
\ben
&& \sum_{m_0+m_1> 0} \prod_{j=0}^{w_1-1} (m_0+w_1m_1-j) \cdot
\prod_{i=2}^n \prod_{j=0}^{w_i-1} (w_i m_1 -j) \tilde{a}_{m_0,m_1}^{(0)} \tilde{x}_0^{m_0} \tilde{x}_1^{m_1} \\
& + & \sum_{m=1}^\infty \sum_{j=0}^{w_1-1} \frac{1}{w_1m-j} \cdot
\frac{(km)!}{\prod_{i=1}^n (w_i(m-1))!} \tilde{x}_1^m \\
& = & \tilde{x}_1 \sum_{m_0+m_1 > 0} \prod_{j=1}^{k} (m_0+km_1+j) \tilde{a}_{m_0,m_1}^{(0)} \tilde{x}_0^{m_0} \tilde{x}_1^{m_1} \\
& + & \tilde{x}_1 \sum_{m=0}^\infty \sum_{j=1}^{k}\frac{1}{km+j}
\frac{(k(m+1))!}{\prod_{i=1}^n (w_im)!} \tilde{x}_1^m.
\een
Compare the coefficients of $\tilde{x}_1$ on both sides:
$$
0 = k! \sum_{j=1}^k \frac{1}{j},
$$
a contradiction.

\section{Integrality of Local Open-closed Mirror Maps}

\subsection{Integrality of local open-closed mirror map: The outer brane case}
Consider the charge vectors
\be \label{eqn:Charge}
\begin{pmatrix}
\hat{l}^{(1)} \\ \hat{l}^{(0)}
\end{pmatrix}
= \begin{pmatrix}
\tilde{l}^{(1)} \\ -\tilde{l}^{(0)}
\end{pmatrix}
= \biggl( \begin{array}{c|cccc|cc}
-k & w_1 & w_2 & \cdots & w_n & 0 & 0 \\
1 &  -1 & 0 & \cdots & 0 & 1 & -1
\end{array} \biggr).
\ee
The associated variables are
\bea
\hat{x}_0 = \frac{a_0a_{n+1}}{a_1a_{n+2}} = \frac{1}{\tilde{x}_0}, \qquad
\hat{x}_1 = (-1)^k \frac{a_1^{w_1} \cdots a_n^{w_n}}{a_0^k} = \tilde{x}_1.
\eea
Suppose that we have a function $\Pi(\hat{x}_0, \hat{x}_1)$ such that
the following equations are satisfied:
\be
\cD_{\hat{l}^{(0)}} \Pi(\hat{x}_0, \hat{x}_1) =0,  \qquad
\cD_{\hat{l}^{(1)}} \Pi(\hat{x}_0, \hat{x}_1) = 0.
\ee
They give rise to the following equations:
\bea
&& \hat{\theta}_0 (\hat{\theta}_0 - k \hat{\theta}_1) \Pi
= \hat{x}_0 \hat{\theta}_0 (\hat{\theta}_0 - w_1 \hat{\theta}_1) \Pi, \\
&& \prod_{j=0}^{w_1-1} (-\hat{\theta}_0 + w_1\hat{\theta}_1 - j) \cdot
\prod_{i=2}^n \prod_{j=0}^{w_i-1} (w_i\hat{\theta}_1 - j) \Pi
= (-1)^k \hat{x}_1 \prod_{j=0}^{k-1} (\hat{\theta}_0 - k \hat{\theta}_1 - j) \Pi.
\eea
We will refer to this system as the local extended Picard-Fuchs system associated with the charge
vectors $\hat{l}^{(0)}$ and $\hat{l}^{(1)}$.

Clearly $\Phi_0 = 1$ is a solution,
and we have the following logarithmic solutions:
\be
\Phi_1^{(0)} = \log \hat{x}_0 + \sum_{m\geq 1} \frac{(km)!}{\prod_{i=1}^n (w_im)!} \frac{\hat{x}_1^m}{m},
\quad \Phi_1^{(1)} = \log \hat{x}_1 - \sum_{m\geq 1} \frac{(km)!}{\prod_{i=1}^n (w_im)!} \frac{\hat{x}_1^m}{km}.
\ee
The local open-closed mirror map is given by
\be
\hat{X}_0 = \hat{x}_0 \exp \sum_{m\geq 1} \frac{(km)!}{\prod_{i=1}^n (w_im)!} \frac{\hat{x}_1^m}{m},
\qquad \hat{X}_1 = \hat{x}_1 \exp \sum_{m\geq 1} \frac{-(km)!}{\prod_{i=1}^n (w_im)!} \frac{\hat{x}_1^m}{km}.
\ee
By a result proved in \cite{Zho1},
we have
\ben
&& \hat{x}_0^{-1}\hat{X}_0, \hat{x}_1^{-1}\hat{X}_1 \in \bZ[[\hat{x}_1]].
\een

\subsection{Integrality of local open-closed mirror map: An inner brane case}
One can also consider the following charge vectors:
\bea
&& \begin{pmatrix} \tilde{\hat{l}}^{(1)} \\ \tilde{\hat{l}}^{(0)} \end{pmatrix}
= \begin{pmatrix} \tilde{l}^{(1)} + \tilde{l}^{(0)} \\ \tilde{l}^{(0)} \end{pmatrix}
= \biggl( \begin{array}{c|cccc|cc}
-k+1&  w_1 -1 & w_2 & \cdots & w_n & -1 & 1 \\
-1 & 1 & 0 & \cdots & 0 & 1 & -1
\end{array} \biggr).
\eea
The associated moduli parameters are
\be
\tilde{\hat{x}}_0 = \frac{a_1a_{n+1}}{a_0a_{n+2}}, \qquad
\tilde{\hat{x}}_1 = (-1)^{k-1} \frac{a_1^{w_1-1}a_2^{w_2} \cdots a_n^{w_n}a_{n+2}}{a_0^{k-1}a_{n+1}}.
\ee
The local extended Picard-Fuchs system associated with these charge vectors is given by:
\bea
&& (\tilde{\hat{\theta}}_0-\tilde{\hat{\theta}}_1)(\tilde{\hat{\theta}}_0 + (w_1-1) \tilde{\hat{\theta}}_1) \Pi
= \tilde{\hat{y}}_0 (-\tilde{\hat{\theta}}_0 + \tilde{\hat{\theta}}_1)
(-\tilde{\hat{\theta}}_0 - (k-1) \tilde{\hat{\theta}}_1) \Pi, \\
&& (\tilde{\hat{\theta}}_0 - \tilde{\hat{\theta}}_1) \cdot
\prod_{j=0}^{w_1-2} (\tilde{\hat{\theta}}_0 + (w_1-1) \tilde{\hat{\theta}}_1-j) \cdot
\prod_{i=2}^n \prod_{j=0}^{w_i-1} (w_i \tilde{\hat{\theta}}_1 - j) \Pi  \\
&& = \hat{y}_1 \cdot  (-\tilde{\hat{\theta}}_0 + \tilde{\hat{\theta}}_1)
\prod_{j=0}^{k-2} (\tilde{\hat{\theta}}_0 + (k-1)\tilde{\hat{\theta}}_1+j) \Pi. \nonumber
\eea
Clearly $\hat{\Phi}_0 = 1$ is a solution and we have two logarithmic solutions:
\bea
&& \hat{\Phi}_1^{(1)}
= \log \tilde{\hat{x}}_1 + \sum_{m=1}^\infty \frac{(km)!}{\prod_{i=1}^n (w_im)!}
\frac{(k-1) (\tilde{\hat{x}}_0 \tilde{\hat{x}}_1)^m}{km}, \\
&&  \hat{\Phi}_1^{(0)}
= \log \tilde{\hat{x}}_0 - \sum_{m=1}^\infty \frac{(km)!}{\prod_{i=1}^n (w_im)!}
\frac{(\tilde{\hat{x}}_0 \tilde{\hat{x}}_1)^m}{km}.
\eea
The local open-closed mirror map is given by
\bea
&& \tilde{\hat{X}}_1 = e^{\hat{\Phi}_1^{(1)}}
= \tilde{\hat{x}}_1 \exp \sum_{m=1}^\infty \frac{(km)!}{\prod_{i=1}^n (w_im)!}
\frac{(k-1) (\tilde{\hat{x}}_0\tilde{\hat{x}}_1)^m}{km}, \\
&& \tilde{\hat{X}}_0 = e^{ \hat{\Phi}_1^{(0) }}
= \tilde{\hat{x}}_0 \exp \sum_{m=1}^\infty \frac{-(km)!}{\prod_{i=1}^n (w_im)!}
\frac{(\tilde{\hat{x}}_0 \tilde{\hat{x}}_1)^m}{km}.
\eea
By a result in \cite{Zho1},
we have
$\tilde{\hat{x}}_0^{-1} \tilde{\hat{X}}_0, \tilde{\hat{x}}_1^{-1} \tilde{\hat{X}}_1 \in \bZ[[\tilde{\hat{x}}_0\tilde{\hat{x}}_1]]$.

\subsection{Integrality of local open-closed mirror map: Another inner brane case}
One can also consider the the charge vectors
\be
\begin{pmatrix}
l^{(1)} \\ l^{(0)}
\end{pmatrix}
= \biggl( \begin{array}{c|cccc|cc}
-k+w_1 & 0 & w_2 & \cdots & w_n & w_1 & -w_1 \\
-1 &  1 & 0 & \cdots & 0 & -1 & 1
\end{array} \biggr),
\ee
where $l^{(1)} = \tilde{l}^{(1)} - w_1 \tilde{l}^{(0)}$ and $l^{(0)} = \tilde{l}^{(0)}$.
The associated moduli parameters are
\be
x_0 = \frac{a_1a_{n+2}}{a_0a_{n+1}}, \qquad
x_1 = (-1)^k \frac{a_2^{w_2} \cdots a_n^{w_n}a_{n+1}^{w_1}}{a_0^{k-w_1} a_{n+2}^{w_1}}.
\ee
The associated local extended Picard-Fuchs system is given by:
\be \label{eqn:PF5}
\check{\cL}_1 S(x_0,x_1) =0, \qquad \check{\cL}_0 S(x_0,x_1) = 0, \qquad \check{\cL}_1' S(x_0,x_1) = 0,
\ee
where
\bea
&& \check{\cL}_0 =  \big[ \theta_0 - x_0 (\theta_0 + (k-w_1) \theta_1) \big] (\theta_0 - w_1 \theta_1), \\
&& \check{\cL}_1 = \prod_{i=2}^n \prod_{j=0}^{w_i-1} (w_i \theta_1-j)
\prod_{j=0}^{w_1-1} (-\theta_0+w_1\theta_1-j) \\
&& \qquad - x_1 \prod_{j=0}^{k-w_1-1} (\theta_0+(k-w_1)\theta_1+j)
\prod_{j=0}^{w_1-1} (-\theta_0+w_1\theta_1+j), \nonumber \\
&& \check{\cL}_1'
= \prod_{j=0}^{w_1-1} (\theta_0-j) \cdot \prod_{i=2}^n \prod_{j=0}^{w_i-1} (w_i \theta_1 -j)
- x_0^{w_1}x_1 \prod_{j=0}^{k-1} (\theta_0+(k-w_1)\theta_1+j).
\eea
We have a holomorphic solution $\check{g}_0 = 1$
and two logarithmic solutions:
\bea
&& \check{g}_1^{(0)} = \log x_0 + \sum_{m=1}^\infty \frac{(km-1)!}{\prod_{i=1}^n (w_im)!} (x_0^{w_1}x_1)^m, \\
&& \check{g}_1^{(1)} = \log x_1 + (k-w_1) \sum_{m=1}^\infty  \frac{(km-1)!}{\prod_{i=1}^n (w_im)!} (x_0^{w_1}x_1)^m.
\eea
The local open-closed mirror map is given by
\bea
&& Q_0 = \exp( \check{g}_1^{(0)}/\check{g}_0) = x_0 \exp \sum_{m=1}^\infty \frac{(km-1)!}{\prod_{i=1}^n (w_im)!} (x_0^{w_1}x_1)^m, \\
&& Q_1 = \exp (\check{g}_1^{(1)}/\check{g}_0) = x_1
\exp \sum_{m=1}^\infty  \frac{(k-w_1)\cdot (km-1)!}{\prod_{i=1}^n (w_im)!} (x_0^{w_1}x_1)^m.
\eea
By a result in \cite{Zho1},
we have
\be
x_0^{-1}Q_0, x_1^{-1}Q_1 \in \bZ[[x_0^{w_1}x_1]].
\ee
Using Lagrange-Good inversion formula as in \cite{Zho1},
one can find explicit expressions for $x_0(Q_0,Q_1), x_1(Q_0,Q_1) \in \bZ[[Q_0,Q_1]]$.
Let $x_i = \sum_{m_0+m_1>0}  A_{m_0,m_1}^{(i)}Q_0^{m_0}Q_1^{m_1}$.
Then $A_{m_0,m_1}^{(0)}$ and $A_{m_0,m_1}^{(1)}$ are the coefficients of
$x_0^{m_0-1}x_1^{m_1}$ and $x_0^{m_0}x_1^{m_1-1}$ respectively of
\ben
\sum_{m=1}^\infty \frac{(km)!}{\prod_{i=1}^n (w_im)!} (x_0^{w_1}x_1)^m
\cdot  \exp (-m_0-(k-w_1)m_1) \sum_{m=1}^\infty \frac{(km-1)!}{\prod_{i=1}^n (w_im)!} (x_0^{w_1}x_1)^m.
\een
It follows that $x_i = Q_i \sum_{a=0}^\infty C_a^{(i)} (Q_0^{w_1}Q_1)^a$,
where $C_a^{(0)}$ is the coefficient of $y^a$ in
\ben
\sum_{m=1}^\infty \frac{(km)!}{\prod_{i=1}^n (w_im)!} y^m
\cdot  \exp (-ka-1) \sum_{m=1}^\infty \frac{(km-1)!}{\prod_{i=1}^n (w_im)!} y^m
\een
and $C_a^{(1)}$ is the coefficient of $y^a$ in
\ben
\sum_{m=1}^\infty \frac{(km)!}{\prod_{i=1}^n (w_im)!} y^m
\cdot  \exp (-ka-(k-w_1)) \sum_{m=1}^\infty \frac{(km-1)!}{\prod_{i=1}^n (w_im)!} y^m.
\een

Now we generalize the discussions in \cite{Zho2} to the open-closed case.
Write
\be
g_0(z_0,z_1) = 1+ \sum_{m_0+m_1>0} c_{m_0,m_1} q_0^{m_0}q_1^{m_1}
=  1+ \sum_{m_0+m_1>0} C_{m_0,m_1} Q_0^{m_0}Q_1^{m_1}.
\ee
If Conjecture 1 holds,
then the coefficients $\{c_{m_0,m_1}\}_{m_0 + m_1\geq 1}$ and $\{C_{m_0,m_1}\}_{m_0+m_1 \geq 1}$  are integers.
One can express $q_0$ and $q_1$ as integral power series
in $Q_0$ and $Q_1$, and vice versa.
If Conjecture 1 holds,
then the following conjecture holds.

\begin{Conjecture}
There are integers $\alpha^{(i)}_{m_0,m_1}$ and $\beta_{m_0,m_1}^{(i)}$ such that
\be
Q_i=q_i\prod_{m_0+m_1\geq 1} (1-q_0^{m_0}q_1^{m_1})^{\alpha^{(i)}_{m_0,m_1}},
\ee
\be
q_i=Q_i\prod_{m_0+m_1\geq 1} (1-Q_0^{m_0}Q_1^{m_1})^{\beta^{(i)}_{m_0,m_1}}.
\end{equation}
for $i=0,1$.
\end{Conjecture}

\section{Examples}

We use Maple to compute some open-closed mirror maps and local open-closed mirror maps and
check their integrality properties for the examples presented below.

\subsection{The $n=2$ case}
There is only one solution $k_1= k_2= 2$ to \eqref{eqn:Sum=1} and so $k=2$ and $w_1=w_2 = 1$.
To the charge vectors
\be \label{eqn:n=2}
\begin{pmatrix}
l^{(1)} \\ l^{(0)}
\end{pmatrix}
= \biggl( \begin{array}{c|cc|cc}
-1 & 0 & 1 & 1 & -1 \\
-1 &  1 & 0 & -1 & 1
\end{array} \biggr)
\ee
the extended Picard-Fuchs operators are:
\ben
&& \cL_0 =  \big[ \theta_0 - x_0 (1 + \theta_0 +  \theta_1) \big] (\theta_0 -  \theta_1), \\
&& \cL_1 =  \theta_1(-\theta_0+w_1\theta_1) - x_1 (\theta_0+\theta_1+1)(-\theta_0+w_1\theta_1),  \\
&& \cL_1' = \theta_0 \theta_1
- x_0 x_1 \prod_{j=1}^{2} (\theta_0+\theta_1+j).
\een
We have the following solutions:
\ben
g_0(x_0,x_1) & = & \sum_{m \geq 0} \frac{(2m)!}{(m!)^2} (x_0x_1)^2, \\
g_1^{(0)}(x_0,x_1) & = & \log x_0 \cdot \sum_{m \geq 0} \frac{(2m)!}{(m!)^2} (x_0x_1)^m
+ \sum_{m=1}^\infty \frac{(2m)!}{(m!)^2} \sum_{j=m+1}^{2m} \frac{1}{j} (x_0x_1)^m \\
& + & \sum_{m_0 \neq m_1}
\frac{(m_0+m_1)!}{m_0! m_1!} \frac{x_0^{m_0} x_1^{m_1}}{m_0-m_1}, \\
g_1^{(1)}(x_0,x_1) & = & \log x_1 \cdot \sum_{m \geq 0} \frac{(2m)!}{(m!)^2} (x_0x_1)^m
+ \sum_{m=1}^\infty \frac{(2m)!}{(m!)^2}
\sum_{j=m+1}^{2m} \frac{1}{j} (x_0x_1)^m  \\
& - & \sum_{m_0 \neq m_1}
\frac{(m_0+m_1)!}{m_0!m_1!} \frac{x_0^{m_0} x_1^{m_1}}{m_0-m_1}.
\een
The first few terms of the Taylor series of the mirror map is given by:
\ben
q_0 & = & x_0 + (x_0^2-x_0x_1) + x_0^3 + (x_0^4 + x_0^3x_1-2x_0^2x_1^2) \\
& + & (x_0^5+2x_0^4x_1-x_0^3x_1^2) +(x_0^6+ 3x_0^5x_1+x_0^4x_1^2-5x_0^3x_1^3) + \cdots, \\
q_1 & = & x_1 + (-x_0x_1+x_1^2)
+x_1^3 + (-2x_0^2x_1^2 + x_0x_1^3+x_1^4) \\
& + & (2x_0x_1^4-x_0^2x_1^3+x_1^5) + (3x_0x_1^5+x_0^2x_1^4-5x_0^3x_1^3+x_1^6) + \cdots.
\een
The inverse open-closed mirror map is given by:
\ben
x_0 & = & q_0 +(- q_0^2 + q_0q_1) + (q_0^3-2q_0^2q_1) + (-q_0^4-3q_0^2q_1^2)  \\
& + & (q_0^5+2q_0^3q_1^2+2q_0^4q_1) + (4q_0^4q_1^2-q_0^6+5q_0^3q_1^3-4q_0^5q_1) + \cdots, \\
x_1 & = & q_1 + (-q_1^2+q_0q_1) + (-2q_0q_1^2+q_1^3) + (-q_1^4-3q_0^2q_1^2) \\
& + & (q_1^5+2q_0^2q_1^3+2q_0q_1^4) + (4q_0^2q_1^4+5q_0^3q_1^3-q_1^6-4q_0q_1^5) + \cdots.
\een

The local extended Picard-Fuchs operators associated with \eqref{eqn:n=2} are:
\ben
&& \check{\cL}_0 =  \big[ \theta_0 - x_0 (\theta_0 + \theta_1) \big] (\theta_0 - \theta_1), \\
&& \check{\cL}_1 = \theta_1 (-\theta_0+\theta_1) - x_1 (\theta_0+\theta_1) (-\theta_0+w_1\theta_1), \\
&& \check{\cL}_1' = \theta_0 \theta_1 - x_0x_1 \prod_{j=0}^{1} (\theta_0+\theta_1+j).
\een
We have a holomorphic solution $\hat{g}_0 = 1$
and two logarithmic solutions:
\ben
&& \check{g}_1^{(0)} = \log x_0 + \sum_{m=1}^\infty \frac{(2m-1)!}{(m!)^2} (x_0x_1)^m
= \log \frac{2x_0}{1+\sqrt{1-4x_0x_1}}, \\
&& \check{g}_1^{(1)} = \log x_1 + \sum_{m=1}^\infty  \frac{(2m-1)!}{(m!)^2} (x_0x_1)^m
= \log \frac{2x_1}{1+\sqrt{1-4x_0x_1}}.
\een
The local open-closed mirror map is given by
\ben
&& Q_0 = x_0 \exp \sum_{m=1}^\infty \frac{(2m-1)!}{(m!)^2} (x_0x_1)^m
= \frac{2x_0}{1 + \sqrt{1-4x_0x_1}}, \\
&& Q_1 = x_1 \exp \sum_{m=1}^\infty  \frac{(2m-1)!}{(m!)^2} (x_0x_1)^m
= \frac{2x_1}{1 + \sqrt{1-4x_0x_1}}.
\een
One can solve these equations to get:
\ben
x_0 = \frac{Q_0}{1+Q_0Q_1}, \qquad
x_1 = \frac{Q_1}{1-Q_0Q_1}.
\een
We have checked
\ben
q_0 = Q_0 \frac{1-Q_1}{1-Q_0}, \qquad q_1 = Q_1 \frac{1-Q_0}{1-Q_1},
\een
and
\ben
Q_0 = q_0 \frac{1+q_1}{1+q_0}, \qquad Q_1 = q_1 \frac{1+q_0}{1+q_1}.
\een

\subsection{The $n=3$ case}

There are three cases to consider:
$(k|w_1, w_2,w_3) = (3|1,1,1)$, $(4|2,1,1)$, $(6|3,2,1)$.
For the first case,
consider the charge vectors:
\ben
\begin{pmatrix}
l^{(1)} \\ l^{(0)}
\end{pmatrix}
= \biggl( \begin{array}{c|ccc|cc}
-2 & 0 & 1 & 1 & 1 & -1 \\
-1 & 1 & 0 & 0 & -1 & 1
\end{array} \biggr),
\een
The extended Picard-Fuchs operators are given by:
\ben
&& \cL_0 =  \big[ \theta_0 - x_0 (1 + \theta_0 + 2 \theta_1) \big] (\theta_0 - \theta_1), \\
&& \cL_1 = \theta_1^2 (-\theta_0+\theta_1)
- x_1 \prod_{j=1}^{2} (\theta_0+2\theta_1+j) (-\theta_0+\theta_1),  \\
&& \cL_1'
=  \theta_0 \theta_1^2
- x_0 x_1 \prod_{j=1}^{3} (\theta_0+2\theta_1+j),
\een
with solutions
\ben
g_0(x_0,x_1) & = & \sum_{m \geq 0} \frac{(3m)!}{(m!)^3} (x_0x_1)^m, \\
g_1^{(0)}(x_0,x_1) & = & \log x_0 \cdot \sum_{m \geq 0} \frac{(3m)!}{(m!)^3} (x_0x_1)^m
+ \sum_{m=1}^\infty \frac{(3m)!}{(m!)^3} \sum_{j=m+1}^{3m} \frac{1}{j}  (x_0x_1)^m \\
& + & \sum_{m_0 \neq m_1}
\frac{(m_0+2m)!}{m_0!(m_1!)^2} \frac{x_0^{m_0} x_1^{m_1}}{m_0-m_1}, \\
g_1^{(1)}(x_0,x_1) & = & \log x_1 \cdot \sum_{m \geq 0} \frac{(3m)!}{(m!)^3} (x_0x_1)^m
+ 2 \sum_{m=1}^\infty \frac{(3m)!}{(m!)^3} \sum_{j=m+1}^{3m} \frac{1}{j}(x_0x_1)^m  \\
& - & \sum_{m_0 \neq m_1}
w_1 \frac{(m_0+2m_1)!}{m_0!(m!)^2} \frac{x_0^{m_0} x_1^{m_1}}{m_0-m_1}.
\een
The open-closed mirror map is given by:
\ben
q_0 & = & x_0 +(x_0^2-2x_0x_1) + (x_0^3+3x_0^2x_1-x_0x_1^2) +
(x_0^4+9x_0^3x_1-29x_0^2x_1^2-2x_0x_1^3) \\
& + & (x_0^5+16x_0^4x_1+27x_0^3x_1^2-23x_0^2x_1^3-5x_0x_1^4) \\
& + & (x_0^6+24x_0^5x_1+127x_0^4x_1^2-527x_0^3x_1^3-63x_0^2x_1^4-14x_0x_1^5) + \cdots, \\
q_1 & = & x_1 + (-x_0x_1 + 2x_1^2)
+ (8x_0x_1^2+5x_1^3)
+ (-16x_0^2x_1^2+33x_0x_1^3 + 14x_1^4) \\
& + & (-x_0^3x_1^2+111x_0^2x_1^3+124x_0x_1^4+42x_1^5) \\
& + & ( -x_0^4x_1^2-307x_0^3x_1^3 +606x_0^2x_1^4 +462x_0x_1^5 + 132 x_1^6) + \cdots.
\een
The local open-closed mirror map is given by:
\ben
&& Q_0 = x_0 \exp \sum_{m=1}^\infty \frac{(3m-1)!}{(m!)^3} (x_0x_1)^m, \\
&& Q_1 = x_1 \exp \sum_{m=1}^\infty  \frac{2 \cdot (3m-1)!}{(m!)^3} (x_0x_1)^m.
\een

For $(k|w_1, w_2,w_3) = (4|2,1,1)$,
one can consider the charge vectors:
\ben
\begin{pmatrix}
l^{(1)} \\ l^{(0)}
\end{pmatrix}
= \biggl( \begin{array}{c|ccc|cc}
-2 & 0 & 1 & 1 & 2 & -2 \\
-1 & 1 & 0 & 0 & -1 & 1
\end{array} \biggr),
\een
The open-closed mirror map is given by:
\ben q_0 & = & x_0 +(x_0^2-x_0x_1) + (x_0^3-7x_0^2x_1-x_0x_1^2) +
(x_0^4+12x_0^3x_1-5x_0^2x_1^2-2x_0x_1^3) \\
& + & (x_0^5+20x_0^4x_1-39x_0^3x_1^2-14x_0^2x_1^3-5x_0x_1^4) \\
& + & (x_0^6+29x_0^5x_1-299x_0^4x_1^2-80x_0^3x_1^3-45x_0^2x_1^4-14x_0x_1^5) + \cdots, \\
q_1 & = & x_1 + (-2x_0x_1 + 2x_1^2) + (x_0^2x_1+8x_0x_1^2+5x_1^3)
+ (4 x_0^2x_1^2+34x_0x_1^3 + 14x_1^4) \\
& + & (- 56x_0^3x_1^2+107x_0^2x_1^3+128x_0x_1^4+42x_1^5) \\
& + & (40 x_0^4x_1^2 +220x_0^3x_1^3 + 592 x_0^2x_1^4 + 476x_0x_1^5 + 132 x_1^6) + \cdots. \een
The inverse open-closed mirror map is given by:
\ben
x_0 & = & q_0 + (-q_0^2+q_0q_1) + (q_0^3+6q_0^2q_1) + (-q_0^2q_1^2-q_0^4-25q_0^3q_1) + \cdots, \\
x_1 & = & q_1 + (2q_0q_1-2q_1^2) + (q_0^2q_1-18q_0q_1^2+3q_1^3) + (50q_0q_1^3-56q_0^2q_1^2-4q_1^4) + \cdots.
\een
The local open-closed mirror map is given by:
\ben
&& Q_0 = x_0 \exp \sum_{m=1}^\infty \frac{(4m-1)!}{(2m)!(m!)^2} (x_0^2x_1)^m, \\
&& Q_1 = x_1 \exp \sum_{m=1}^\infty  \frac{2\cdot (4m-1)!}{(2m)!(m!)^2} (x_0^2x_1)^m.
\een
The inverse local open-closed mirror map is given by:
\ben
&& x_0 = Q_0(1-3Q_0^2Q_1 - 12 Q_0^4Q_1^2 - 253 Q_0^6Q_1^3 - 6033 Q_0^8Q_1^4+ \cdots), \\
&& x_1 = Q_1(1-6Q_0^2Q_1 -15 Q_0^4Q_1^2 -434 Q_0^6Q_1^3 -10404 Q_0^8Q_1^4+ \cdots).
\een
One has
\ben
q_0 & = & Q_0-Q_0Q_1+Q_0^2-Q_0Q_1^2-7Q_0^2Q_1+Q_0^3-2Q_0Q_1^3-5Q_0^2Q_1^2+9Q_0^3Q_1+Q_0^4\\
& & -5Q_0Q_1^4 -14 Q_0^2Q_1^3 -30 Q_0^3Q_1^2 + 14 Q_0^4Q_1 + Q_0^5 \\
& & - 14Q_0Q_1^5-45Q_0^2Q_1^4
-65Q_0^3Q_1^3-215Q_0^4Q_1^2+20Q_0^5Q_1+Q_0^6 + \cdots \\
q_1 & = & Q_1+2Q_1^2-2Q_0Q_1+5Q_1^3+8Q_0Q_1^2+Q_0^2Q_1+14Q_1^4+34Q_0Q_1^3-2Q_0^2Q_1^2\\
&& +42Q_1^5+128Q_0Q_1^4+83Q_0^2Q_1^3-38Q_0^3Q_1^2\\
&& +132Q_1^6+476Q_0Q_1^5+502Q_0^2Q_1^4+100Q_0^3Q_1^3
+28Q_0^4Q_1^2+1776Q_0Q_1^6 + \cdots.
\een
One can also consider the charge vectors corresponding to another phase of the $B$-brane:
\ben
\begin{pmatrix}
l^{(1)} \\ l^{(0)}
\end{pmatrix}
= \biggl( \begin{array}{c|ccc|cc}
-3 & 0 & 2 & 1 & 1 & -1 \\
-1 & 1 & 0 & 0 & -1 & 1
\end{array} \biggr),
\een
Then the open-closed mirror map is given by:
\ben
q_0 & = & x_0 + (x^2-3x_0x_1) + (x_0^3+10x_0^2x_1-3x_0x_1^2) \\
& + & (x^4+28x_0^3x_1-111x_0^2x_1^2-10x_0x_1^3) + \cdots, \\
q_1 & = & x_1 + (-x_0x_1+3x_1^2) + (42x_0x_1^2+12x_1^3) \\
& + & (-63x_0^2x_1^2+192x_0x_1^3+55x_1^4) + \cdots.
\een

For $(k|w_1, w_2,w_3) = (6|3,2,1)$,
we check the following three cases corresponding to three different phases of the $B$-brane:
For the charge vectors
\ben
\begin{pmatrix}
l^{(1)} \\ l^{(0)}
\end{pmatrix}
= \biggl( \begin{array}{c|ccc|cc}
-3 & 0 & 2 & 1 & 3 & -3 \\
-1 & 1 & 0 & 0 & -1 & 1
\end{array} \biggr),
\een
we have
\ben
q_0 & = & x_0 + (x_0^2-x_0x_1) + (x_0^3-7x_0^2x_1-2x_0x_1) \\
& + & (x_0^4-37x_0^3x_1-17x_0^2x_1^2-7x_0x_1^3) + \cdots, \\
q_1 & = & x_1 + (-3x_0x_1+3x_1^2) + (3x_0^2x_1+9x_0x_1^2+12x_1^3) \\
& + & (-x_0^3x_1+45x_0^2x_1^2+81x_0x_1^3+55x_1^4) + \cdots;
\een
for the charge vectors
\ben
\begin{pmatrix}
l^{(1)} \\ l^{(0)}
\end{pmatrix}
= \biggl( \begin{array}{c|ccc|cc}
-4 & 0 & 3 & 1 & 2 & -2 \\
-1 & 1 & 0 & 0 & -1 & 1
\end{array} \biggr),
\een
we have
\ben
q_0 & = & x_0 + (x_0^2-2x_0x_1) + (x_0^3-22x_0^2x_1-5x_0x_1^2) \\
& + & (x_0^4+212x_0^3x_1-49x_0^2x_1^2-24x_0x_1^3) + \cdots, \\
q_1 & = & x_1 + (-2x_0x_1+4x_1^2) + (x_0^2x_1+32x_0x_1^2+22x_1^2) \\
& + & (342x_0^2x_1^2+284x_0x_1^3+140x_1^4) + \cdots;
\een
for the charge vectors
\ben
\begin{pmatrix}
l^{(1)} \\ l^{(0)}
\end{pmatrix}
= \biggl( \begin{array}{c|ccc|cc}
-5 & 0 & 2 & 3 & 1 & -1 \\
-1 & 1 & 0 & 0 & -1 & 1
\end{array} \biggr),
\een
we have
\ben
q_0 & = & x_0 + (x_0^2-10x_0x_1) + (x_0^3+77x_0^2x_1-55x_0x_1^2) \\
& + & (x_0^4+227x_0^3x_1-2635x_0^2x_1^2-785x_0x_1^3) + \cdots, \\
q_1 & = & x_1 + (-x_0x_1 + 10x_1^2)
+(525x_0x_1^2+155x_1^3) \\
& + & (-685x_0^2x_1^2+6905x_0x_1^3+2885x_1^4) + \cdots.
\een

\subsection{Some $n>3$ cases}

For $n=4$,
there are thirteen cases to consider when $w_1 \geq \cdots \geq w_4$ listed in \cite{Zho2};
when $n=5$ there are 147 cases.
We have checked some of them using our Maple algorithm.
The following cases are particularly interesting because they correspond to the large volume phase
of $B$-brane geometries in one-moduli Calabi-Yau hypersurfaces in weighted
projective spaces:
$$(k|w_1, \dots, w_5) = (5|1,1,1,1,1), (6|2,1,1,1,1),
(8|4,1,1,1,1),(10|5,2,1,1,1).
$$
When $(k|w_1, \dots, w_5) = (5|1,1,1,1,1)$, for the charge vectors
\be
\begin{pmatrix} l^{(1)} \\ l^{(0)} \end{pmatrix}
= \biggl( \begin{array}{c|ccccc|cc}
-4 &  0 & 1 & 1 & 1 & 1&  -1 & 1 \\
-1 & 1 & 0 & 0 & 0 & 0 &1 & -1
\end{array} \biggr)
\ee
we have
\ben
q_0 & = & x_0 + (x_0^2-24x_0x_1) +(x_0^3-24x_0^2x_1-972x_0x_1^2) \\
& + & (x_0^4+216x_0^3x_1-20772x_0^2x_1^2-95264x_0x_1^3) + \cdots, \\
q_1 & = & x_1 + (-x_0x_1+24x_1^2) + (746x_0x_1^2+1548x_1^3) \\
& + & (-1010 x_0^2x_1^2 + 36732 x_0x_1^3 + 155744 x_1^4) + \cdots.
\een

When $(k|w_1, \dots, w_5) = (6|2,1,1,1,1)$, we consider two phases.
For the charge vectors
\ben
\begin{pmatrix} l^{(1)} \\ l^{(0)} \end{pmatrix}
= \biggl( \begin{array}{c|ccccc|cc}
-4 & 0 & 1 & 1 & 1 & 1 & -2 & 2 \\
-1 & 1 & 0 & 0 & 0 & 0 & 1 & -1
\end{array} \biggr)
\een
we have
\ben
q_0 & = & x_0 + (x_0^2-12x_0x_1) + (x_0^3-132x_0^2x_1-558x_0x_1^2) \\
& + & (x_0^4+570x_0^3x_1-6678x_0^2x_1^2-54328x_0x_1^3) + \cdots, \\
q_1 & = & x_1 + (-2x_0x_1+24x_1^2) + (x_0^2x_1+192x_0x_1^2+1548x_1^3) \\
& + & (1632x_0^2x_1^2+17784x_0x_1^3+155744x_1^4) + \cdots;
\een
for the charge vectors
\ben
\begin{pmatrix} l^{(1)} \\ l^{(0)} \end{pmatrix}
= \biggl( \begin{array}{c|ccccc|cc}
-5 &  0 & 2 & 1 & 1 & 1&  -1 & 1 \\
-1 & 1 & 0 & 0 & 0 & 0 &1 & -1
\end{array} \biggr)
\een
we have
\ben
q_0 & = & x_0 + (x_0^2-60x_0x_1) + (x_0^3+462x_0^2x_1-7650x_0x_1^2) \\
& + & (x_0^4+1362x_0^3x_1-225270x_0^2x_1^2-2271800x_0x_1^3) + \cdots, \\
q_1 & = & x_1 + (-x_0x_1+60x_1^2)+(-60x_0x_1^2+11250x_1^3) \\
& + & (1890x_0^2x_1^2+175050x_0x_1^3+3405800x_1^4) + \cdots.
\een

When $(k|w_1, \dots, w_5) = (8|4,1,1,1,1)$,
for the charge vectors
\ben
\begin{pmatrix} l^{(1)} \\ l^{(0)} \end{pmatrix}
= \biggl( \begin{array}{c|ccccc|cc}
-4 &  0 & 1 & 1 & 1 & 1&  -4 & 4 \\
-1 & 1 & 0 & 0 & 0 & 0 &1 & -1
\end{array} \biggr)
\een
we have
\ben
q_0 & = & x_0 + (x_0^2-6x_0x_1) + (x_0^3-46x_0^2x_1-297x_0x_1^2) \\
& + & (x_0^4-226x_0^3x_1-3297x_0^2x_1^2-28946x_0x_1^3) + \cdots, \\
q_1 & = & x_1 + (-4x_0x_1+24x_1^2) + (6x_0^2x_1+64x_0x_1^2+1548x_1^3) \\
& + & (-4x_0^3x_1+224 x_0^2x_1^2 + 10608 x_0x_1^3 + 155744 x_1^4) + \cdots;
\een
for the charge vectors
\be
\begin{pmatrix} l^{(1)} \\ l^{(0)} \end{pmatrix}
= \biggl( \begin{array}{c|ccccc|cc}
-7 &  0 & 4 & 1 & 1 & 1 &  -1 & 1 \\
-1 & 1 & 0 & 0 & 0 & 0 & 1 & -1
\end{array} \biggr)
\ee
we have
\ben
q_0 & = & x_0 + (x_0^2-210x_0x_1) + (x_0^3 + 2676x_0^2x_1 -113085x_0x_1^2) \\
& + & (x_0^4+8556x_0^3x_1-4420395x_0^2x_1^2-137765950x_0x_1^3) + \cdots, \\
q_1 & = & x_1 + (-x_0x_1+210x_1^2) + (23212x_0x_1^2+157185x_1^3) \\
& + & (-29302x_0^2x_1^2+8462685x_0x_1^3+194522650x_1^4) + \cdots.
\een

When $(k|w_1, \dots, w_5) = (10|5,2,1,1,1)$,
for the charge vectors
\ben
\begin{pmatrix} l^{(1)} \\ l^{(0)} \end{pmatrix}
= \biggl( \begin{array}{c|ccccc|cc}
-5 &  0 & 2 & 1 & 1 & 1&  -5 & 5 \\
-1 & 1 & 0 & 0 & 0 & 0 &1 & -1
\end{array} \biggr)
\een
we have
\ben
q_0 & = & x_0 + (x_0^2-12x_0x_1) + (x_0^3-102x_0^2x_1-1818x_0x_1^2) \\
& + & (x_0^4-522x_0^3x_1-23838x_0^2x_1^2-538168x_0x_1^3) + \cdots, \\
q_1 & = & x_1 + (-5x_0x_1+60x_1^2) + (10x_0^2x_1 + 150 x_0x_1^2 + 11250 x_1^3) \\
& + & (-10x_0^3x_1 + 450 x_0^2x_1^2 + 86250 x_0x_1^3 + 3405800 x_1^4) + \cdots;
\een
for the charge vectors
\be
\begin{pmatrix} l^{(1)} \\ l^{(0)} \end{pmatrix}
= \biggl( \begin{array}{c|ccccc|cc}
-8 &  0 & 5 & 1 & 1 & 1&  -2 & 2 \\
-1 & 1 & 0 & 0 & 0 & 0 &1 & -1
\end{array} \biggr)
\ee
we have
\ben
q_0 & = & x_0 + (x_0^2-168x_0x_1) + (x_0^3-3192x_0^2x_1-166068x_0x_1^2) \\
& + & (x_0^4+33534x_0^3x_1-3742116x_0^2x_1^2-336621152x_0x_1^3) + \cdots, \\
q_1 & = & x_1 + (-2x_0x_1 + 336x_1^2) + (x_0^2x_1 + 5376 x_0x_1^2 + 416808 x_1^3) \\
& + & (-11760 x_0^2x_1^2 + 9366672 x_0x_1^3 + 859605376 x_1^4) + \cdots;
\een
for the charge vectors
\be
\begin{pmatrix} l^{(1)} \\ l^{(0)} \end{pmatrix}
= \biggl( \begin{array}{c|ccccc|cc}
-9 & 0&  5 & 2  & 1 & 1&  -1 & 1 \\
-1 & 1 & 0 & 0 & 0 & 0 &1 & -1
\end{array} \biggr)
\ee
we have
\ben
q_0 & = & x_0 + (x_0^2-1512x_0x_1) + (x_0^3+ 27654 x_0^2x_1 - 8046108 x_0x_1^2) \\
& + & (x_0^4 + 95694 x_0^3 x_1 - 378472500 x_0^2x-1^2 - 93766645728 x_0 x_1^3) + \cdots, \\
q_1 & = & x_1 + (-x_0x_1 + 1512 x_1^2) + (265332x_0x_1^2 + 10332252 x_1^3) \\
& + & (121554726048x_1^4+719463276x_0x_1^3-334884x_0^2x_1^2) + \cdots.
\een
Some $n>4$ cases have also been checked.

\vspace{.1in}
{\em Acknowledgements}.
This research is partly supported by NSFC grants (10425101 and 10631050)
and a 973 project grant NKBRPC (2006cB805905).

\end{document}